\documentclass[12pt,a4paper,oneside,final,reqno]{amsart}

\usepackage{amsthm}
\usepackage{amsmath}
\usepackage{amssymb}
\usepackage{mathtools}
\usepackage{amscd}
\usepackage[utf8]{inputenc}
\usepackage{typearea}

\usepackage{yfonts}
\usepackage{textcomp}
\usepackage{mathrsfs}
\usepackage{hyperref}
\usepackage[draft]{fixme} 
\usepackage{pdfsync}
\usepackage[active]{srcltx}
\usepackage{verbatim}
\usepackage{todonotes}
\usepackage[capitalise]{cleveref}
\usepackage{tikz-cd}
\usepackage{tensor}

\usepackage[backend=bibtex, style=alphabetic, maxnames=50]{biblatex}

\newcommand{\N}{\mathbb{N}}

\newcommand{\Go}{G^{(0)}}

\newcommand{\dad}{\operatorname{dad}}

\theoremstyle{plain}
\newtheorem{theorem}{Theorem}[section]
\newtheorem{corollary}[theorem]{Corollary}
\newtheorem{lemma}[theorem]{Lemma}

\newtheorem{introtheorem}{Theorem}

\newtheorem{introcorollary}[introtheorem]{Corollary}

\theoremstyle{definition}
\newtheorem{definition}[theorem]{Definition}

\newtheorem{remark}[theorem]{Remark}

\numberwithin{equation}{section}
 
\title{Polynomial growth and dynamic asymptotic dimension}

\date{}

\author[1]{Are Austad}
\address{Are Austad, Department of Mathematics, University of Oslo, P.O. Box 1053 Blindern, 0316 Oslo, Norway}
\email{areaus@math.uio.no}

\author[2]{Christian Bönicke}
\address{Christian Bönicke, School of Mathematics, Statistics and Physics, Newcastle University, Newcastle upon Tyne NE1 7RU,
UK}
\email{christian.bonicke@ncl.ac.uk}


\begin{document}
\begin{abstract}
   We show that principal ample groupoids with local polynomial growth of degree at most $s$ have dynamic asymptotic dimension at most $\lfloor s\rfloor$. As a consequence, $\sigma$-compact principal ample groupoids with local polynomial growth have slow dynamic asymptotic dimension growth.
\end{abstract}

\maketitle
\section{Introduction}
The dynamic asymptotic dimension is a measure of the dynamical complexity of a topological dynamical system. It was introduced in \cite{GWY17} for actions of discrete groups on locally compact spaces, and more generally for locally compact Hausdorff \'etale groupoids to produce an upper bound for the nuclear dimension of the associated groupoid $C^*$-algebra (see also \cite{CDGaHV24,BL23} for extensions of this result). These upper bounds for the nuclear dimension are often far from optimal \cite{CETWW21}. However, in \cite{LLW26} the authors introduce a version of nuclear dimension for an inclusion $D\subseteq A$ of $C^*$-algebras, called diagonal dimension, which seems to be very closely related to the dynamic asymptotic dimension of the groupoid model for the diagonal inclusion $D\subseteq A$.
In a different direction it was shown in \cite{BDGW23} that the dynamic asymptotic dimension has a homological shadow: 
for a $\sigma$-compact principal ample groupoid $G$, its homology vanishes in degrees above $\dad(G)$.

The main contribution of the present article is a concrete sufficient criterion for finite dynamic asymptotic dimension (resp. control over the dynamic asymptotic dimension growth) in terms of the local polynomial growth. Roughly, local polynomial growth means that, for each compact open symmetric set $D$ of arrows, the number of arrows in $D^n$ with any fixed source is bounded by a polynomial in $n$, uniformly over the source. The polynomial may depend on $D$. Degree at most $s$ means that its exponent can always be chosen to be $s$. 
\begin{introtheorem}\label{introthm:A}
    Let $G$ be a principal ample groupoid with local polynomial growth of degree at most $s\geq 0$. Then $\dad(G)\leq \lfloor s\rfloor$. 
\end{introtheorem}

The present work is inspired by \cite{P23} where a similar statement is proved in the context of graphs and \cite{BY25} who prove a similar statement in the Borel setting. 
In fact, 
via the coarse groupoid construction, Theorem \ref{introthm:A} recovers Papasoglu's bound $\operatorname{asdim}(X)\leq k$ when the uniform growth of $X$ is $O(r^s)$ for some fixed $s<k+1$, whereas \cite{P23} requires only the weaker hypothesis $o(r^{k+1})$. 
Moreover, the Borel asymptotic dimension result \cite[Theorem~1.25]{BY25} can be obtained by following our strategy as well, making minor modifications to the arguments presented below.

Finally, in \cite{WWZZ26}, the authors introduce the concept of dynamic asymptotic dimension growth as a tool to study groupoids with infinite dynamic asymptotic dimension.
\begin{introcorollary}
    Let $G$ be a $\sigma$-compact, principal, ample groupoid with local polynomial growth and $f:[0,\infty)\to [0,\infty)$ an increasing unbounded function. Then there exists a continuous coarse length function $\ell$ on $G$ such that $\dad_{G,\ell}(t)\leq f(t)$. In particular, $\ell$ can be chosen so that $(G,\ell)$ has slow dynamic asymptotic dimension growth. 
    \end{introcorollary}
In Section \ref{sec:examples} we apply Theorem \ref{introthm:A} to relevant examples. Beyond recovering known dynamic asymptotic dimension bounds in the case of coarse groupoids and AF groupoids, we also consider partial action groupoids and point set groupoids for point sets in locally compact groups. Specifically, we obtain the bound $\lfloor s \rfloor$ for free partial actions of discrete groups with polynomial growth of degree at most $s$ on locally compact totally disconnected spaces. We also establish the same bound for point set groupoids, assuming separation, finite local complexity, trivial translation stabilizers throughout the non-empty hull, and polynomial volume growth of degree at most $s$ in the ambient group.

\subsubsection*{Acknowledgments} 
The first author wishes to thank Newcastle University for its hospitality during a research stay in February 2026, during which parts of this research were conducted. 
The first author was funded by the Research Council of Norway projects 324944 and 359602.

\section{Groupoids and dynamic asymptotic dimension}\label{sec: prelims}
Let $G$ denote a groupoid with unit space $\Go$, and let $G^{(2)}$ denote its composable pairs. It comes equipped with range and source maps $r,s \colon G \to \Go$. The range and source fibres will be denoted by $G^x = r^{-1}(x)$ and $G_x = s^{-1}(x)$ for $x \in \Go$.
For a subset $V \subseteq \Go$, the restriction of $G$ to $V$ is $G\vert_V := s^{-1}(V) \cap r^{-1}(V)$. We say that $G$ is \emph{principal} if $G\vert_{\{x\}} = s^{-1}(\{x\}) \cap r^{-1}(\{x\}) = \{x \}$ for all $x \in \Go$. 
Furthermore, for subsets $A, B \subseteq G$, we denote the set of all possible products of elements in $A$ and elements in $B$ by $AB = \{g h \in G \mid g \in A, h \in B, (g, h) \in G^{(2)}\}$. All groupoids considered in this article will be  locally compact and Hausdorff. That is, we equip $G$ with a locally compact Hausdorff topology such that, when $\Go$ is equipped with the relative topology from $G$, and $G^{(2)}$ is equipped with the relative topology from $G \times G$, the multiplication and inversion maps on $G$ are continuous. By continuity of composition, if $A,B \subseteq G$ are compact, then $AB$ is also compact. For $A \subseteq G$ we let $\langle A \rangle \subseteq G$ be the algebraically generated subgroupoid generated by $A$.

An \emph{open bisection} is an open subset $B \subseteq G$ on which both $s$ and $r$ restrict to homeomorphisms onto open subsets of $\Go$. The groupoid $G$ is said to be \'etale if $s$ (equivalently $r$) is a local homeomorphism, and is said to be \emph{ample} if its topology has a basis consisting of compact open bisections. 
In an \'etale groupoid, products of open sets are open, and so in particular, finite products of compact open sets are again compact open. Since we exclusively work in Hausdorff topologies, such sets are also closed.

\begin{definition}
    Let $G$ be an \'etale groupoid. We say that $G$ has \emph{local polynomial growth} if for every compact open symmetric set $D\subseteq G$, there exists an $s\geq 0$ and a constant $C_D\geq 1$ such that
    $$\sup_{x\in G^{(0)}} \vert D^n\cap G_x\vert \leq C_D(1+n)^s \quad \text{ for all } n\in \mathbb N_0.$$
    If the same exponent $s\geq 0$ can be chosen independently of $D$, then we say that $G$ has \emph{local polynomial growth of degree at most $s$}.
\end{definition}

\begin{remark} \label{rem:length-local-growth}
We note here that if $G$ admits a controlled length function $\ell:G\to [0,\infty)$ with polynomial growth of degree $s\geq 0$ in the sense of \cite{AustadBonicke26}, then $G$ has local polynomial growth of degree at most $s$. Indeed, write $B_\ell(t) := \{g \in G \mid \ell(g) \leq t \}$ for the closed balls for $\ell$. 
 Let $C \geq 1$ be such that
 \begin{align*}
        \sup_{x \in \Go} \vert B_\ell(t) \cap G_x \vert \leq C (1+t) ^s \qquad \text{for all $t \geq 0$, }
    \end{align*}
    and let $D$ be a compact open symmetric subset of $G$. Define $M_D := \max \{1 , \sup_{g\in D} \ell(g)\}$, which is finite as $\ell$ is assumed controlled. Since $\ell$ is a length function, we have $D^n \subseteq B_\ell(n M_D)$. Then
    \begin{align*}
        \sup_{x \in\Go} \vert D^n \cap G_x \vert \leq \sup_{x \in \Go} \vert  B_\ell(n M_D) \cap G_x \vert \leq C (1 +nM_D)^s \leq CM_D^s (1+n)^s,
    \end{align*}
\end{remark}

Lastly in this section we recall the definition of dynamic
asymptotic dimension.

\begin{definition}\label{def:dad}
    Let $G$ be an \'etale groupoid. We say $G$ has \emph{dynamic asymptotic dimension} bounded by $d$, $d \in \N_0$, and write $\dad(G) \leq d$, if for every relatively compact open set $E \subseteq G$, there are open subsets $U_0 , \ldots , U_d \subseteq \Go$ covering $s(E) \cup r(E)$ such that $\langle E \cap G\vert_{U_i} \rangle$ is relatively compact in $G$ for each $i$. The least such $d$, or $\infty$ if none exists, is the dynamic asymptotic dimension $\dad(G)$ of $G$.
\end{definition}

Our arguments in Section \ref{sec:main-results} will be simplified in the case of a compact unit space. Given an \'etale groupoid $G$ we may pass to its Alexandrov groupoid \cite[Definition 7.8]{KKLRU2021}. Suppose $\Go$ is noncompact. Let $(\Go)^+ = \Go \sqcup \{\infty\}$ be the one-point compactification of the unit space. The Alexandrov groupoid $G^+$ is the set $G \cup (\Go)^+$, with topology determined by requiring that the inclusions $G, (\Go)^+ \subseteq G^+$ are open, and the groupoid structure extends that of $G$ by making $\infty$ a unit. Since $G$ is locally compact Hausdorff \'etale, so is $G^+$. Moreover, if $G$ is ample, so is $G^+$.
We remark that principality and local polynomial growth pass to $G^+$, with the same exponent whenever a uniform exponent is available. 

\begin{lemma}[{\cite[Proposition 3.13]{CDGaHV24}}]\label{lemma:cpct-unit-space-suffices}
    Suppose that $G$ is an étale groupoid. Then $\dad (G) = \dad(G^+)$.
\end{lemma}

When $G$ is an \'etale groupoid with compact unit space and $G$ is equipped with a coarse continuous length function $\ell:G\to [0,\infty)$, then we follow \cite{WWZZ26} and define the dynamic asymptotic dimension function
$$\mathrm{dad}_{G,\ell}:[0,\infty)\to \mathbb N_0 \cup \{\infty \}$$
by letting $\mathrm{dad}_{G,\ell}(t)$ be the least $d\in \mathbb N_0$ such that there are open subsets $U_0 , \ldots , U_d \subseteq \Go$ covering $s(B_\ell(t))$ such that $\langle B_\ell(t) \cap G\vert_{U_i} \rangle$ is relatively compact in $G$ for each $i$, and setting $\dad_{G,\ell}(t) = \infty$ if no such $d$ exists. The growth type of the function $\mathrm{dad}_{G,\ell}$ is called the dynamic asymptotic dimension
growth of $(G,\ell)$. For noncompact unit space, following \cite[Definition 4.10]{WWZZ26} we define $\dad_{G,\ell} := \dad_{G^+, \ell^+}$, where $\ell^+|_G = \ell$ and $\ell^+(\infty) = 0$. This extension is again continuous and coarse.  

\section{Main Results}
\label{sec:main-results}

Before launching into the details of the proof of Theorem~\ref{introthm:A}, let us outline the main ideas and key steps in the proof.
A key ingredient in the proof is to show that for every $N\in \mathbb N$ and every compact open subset $A\subseteq G^{(0)}$ there exists a compact open set $Y\subseteq A$ such that
\begin{equation}\label{eq:separation}
     \langle D\cap G|_{A\setminus Y}\rangle \subseteq D^{2N}.
\end{equation}
A set $Y$ with the properties above plays the role of the separating subsets in \cite{P23}, hence we will sometimes refer to $Y$ as a \emph{separating set} in $A$. As outlined in \cite[Section~3.3]{BY25} and the references therein, such a separating subset can be produced via an explicit ball carving construction associated with suitable radius choice functions $t:A\to \{0,\ldots, N-1\}$. If $t$ happens to be locally constant, the set $Y$ obtained from the ball carving procedure is compact open and satisfies \eqref{eq:separation}.

The strategy to prove Theorem~\ref{introthm:A} is to successively apply this procedure to obtain a decreasing sequence of compact open subsets $G^{(0)}=A_0\supseteq A_1\supseteq A_2\supseteq \ldots$, each 
$A_{i+1}$ separating in $A_i$,
i.e. $\langle D\cap G|_{A_i\setminus A_{i+1}}\rangle\subseteq D^{2N}$. Thus, the goal is to use local polynomial growth to ensure that this process eventually halts in the sense that $A_i$ is empty at some finite stage. 
This is where our argument significantly differs from both the minimality argument in \cite{P23} and the probabilistic approach taken in \cite{BY25}.
We use the fact that local polynomial growth of $G$ implies that 
suitable balls in the orbit graph with respect to $D$ have small relative boundaries 
to produce 
for every $\varepsilon > 0$ and all sufficiently large $N \in \N$ 
a locally constant radius choice function $t:G^{(0)}\to \{0,\ldots, N\}$ such that in addition to \eqref{eq:separation} the set 
$Y$ 
produced by the ball carving construction also has a small growth function $\gamma_Y$ (see \eqref{eq:def-growth-function-unit-space-subsets} below for the definition) relative to that of $A$, in the sense that
\begin{equation}
    \gamma_Y(R)\leq \varepsilon \gamma_A(R+2N).
\end{equation}
Applying this estimate iteratively to our decreasing sequence of subsets $A_0\supseteq A_1\supseteq A_2\supseteq \ldots$ and using the polynomial growth assumption we conclude that $\gamma_{A_{\lfloor s\rfloor +1}}(0)<1$ and hence $A_{\lfloor s\rfloor +1}=\emptyset$. 

Let us now turn to the details of the proof.
By Lemma \ref{lemma:cpct-unit-space-suffices}, we may assume that the unit space is compact, and so the auxiliary lemmas in the sequel will be stated and proved under that assumption. The following remark tells us that in the definition of dynamical asymptotic dimension in Definition \ref{def:dad} it suffices to consider precompact sets which contain the compact unit space.  
\begin{remark}\label{remark:D-containing-Go-suffices}
    When $\Go$ is compact, the definition of dynamical asymptotic dimension may be tested on symmetric relatively compact open sets containing $\Go$, see \cite[Definition 2.9]{Bonicke24} and the surrounding discussion. If $G$ is ample, it moreover suffices to restrict to compact open sets. Indeed, given a relatively
compact open $E\subseteq G$, choose compact open bisections
$B_1,\ldots,B_n$ covering $\overline E$ and set
\begin{align*}
 D=\Go \cup\bigcup_{j=1}^{n}(B_j\cup B_j^{-1}).
\end{align*}
Then $D$ is compact open, symmetric, and contains both $\Go$ and $E$. Any DAD cover $U_0, \ldots , U_d$ for $D$ works for $E$, as $\langle E \cap G\vert_{U_i}\rangle \subseteq \langle D \cap G\vert_{U_i} \rangle$ for all $i$. 

Thus to prove $\dad (G) \leq d$, it suffices to construct, for every such $D$, an open cover $U_0, \ldots, U_d$ of $\Go$ for which $\langle D \cap G\vert_{U_i} \rangle$ is relatively compact for every $0 \leq i \leq d$. 
\end{remark}

In light of Remark \ref{remark:D-containing-Go-suffices}, 
we will focus on compact open symmetric sets $D \subseteq G$ with $D \supseteq \Go$. Indeed, fix such $D$, set $D^0 = \Go$, and define
\begin{align*}
    B_D (x,n) = r(D^n \cap G_x). 
\end{align*}
We identify $G^{(0)}$ with the vertices of a graph, where two units $x,y$ are connected by an edge exactly when there is $g \in D$ with $s(g) = x$ and $r(g) = y$. We then naturally have an extended graph metric on $\Go$, and since $\Go \subseteq D = D^{-1}$, the set $B_D(x,n)$ is exactly the radius $n$ ball about $x\in \Go$ in this extended metric. The distance between different graph components is set to be infinite. 

Moreover, for any $A \subseteq \Go$ we define
\begin{align}\label{eq:def-growth-function-unit-space-subsets}
    \gamma_A (n) = \sup_{x \in \Go} \vert A \cap B_D(x,n)\vert ,
\end{align}
where we remark that the supremum is taken over $\Go$, not just $A$. A \emph{$D$-component of $A$} is a connected component of the induced graph on $A$.

Note the trivial observation
\begin{align*}
    \gamma_A (0) = \begin{cases}
        1 & A \neq \emptyset \\
        0 & A = \emptyset. 
    \end{cases}
\end{align*}
If $G$ is principal, the range map is injective on $G_x$ for every $x \in \Go$, and so
\begin{align*}
    \vert B_D (x,n) \vert = \vert D^n \cap G_x \vert .
\end{align*}

\begin{lemma}\label{lemma:locally-constant-counts}
    Let $G$ be a principal ample groupoid with compact unit space $\Go$. 
    If $T \subseteq G$ and $W \subseteq \Go$ are compact open, then $T \cap r^{-1}(W)$ and $r(T \cap s^{-1}(W))$ are compact open. Moreover, the function
    \begin{align*}
        x \mapsto \vert T \cap G_x \vert 
    \end{align*}
    is bounded and locally constant. In particular, the map $x \mapsto \vert W \cap B_D(x,n) \vert$ is locally constant for each fixed $n$. 
\end{lemma}

\begin{proof}
    Since $r$ and $s$ are local homeomorphisms, they are both continuous and open maps. The first claim is then immediate from $T$ and $W$ being compact open. 

    Using compactness of $T$, cover $T$ by finitely many compact open bisections. Using the fact that $G$ is ample, find a finite disjoint refinement, i.e., find finitely many disjoint compact open bisections $T_1, \ldots , T_k$. Then we may write
    \begin{align*}
        x \mapsto \vert T \cap G_x \vert = \sum_{j=1}^k 1_{s(T_j)} (x)
    \end{align*}
    which is locally constant. The last assertion then follows by setting $T = D^n \cap r^{-1}(W)$ and calculating
    \begin{align*}
        \vert r^{-1}(W) \cap D^n \cap G_x \vert = \vert r(r^{-1}(W) \cap D^n \cap G_x ) \vert  = \vert W \cap B_D(x,n) \vert 
    \end{align*}
    where the first equality follows by injectivity of $r$ on each source fibre. 
\end{proof}

\begin{lemma}\label{lemma:loc-to-glob-D-boundedness}
    Let $G$ be a principal ample groupoid with compact unit space $\Go$, and let $A \subseteq \Go$. If every $D$-component of $A$ is contained in some $B_D(u,N)$, where $u$ is allowed to depend on the component, then
    \begin{align*}
        \langle D\cap G|_A\rangle \subseteq D^{2N}. 
    \end{align*}
\end{lemma}

\begin{proof}
    If $h \in \langle D \cap G|_A \rangle$, then $s(h)$ and $r(h)$ are in the same $D$-component of $A$. This component is contained in $B_D(u,N)$ for some $u$, so choose $p, q \in D^N$ with $s(p) = s(q) = u$ and $r(p) = s(h)$ and $r(q) = r(h)$. Then $s(qp^{-1}) = s(h)$ and $r(qp^{-1}) = r(h)$, so principality yields $h = qp^{-1} \in D^{2N}$. 
\end{proof}

\begin{lemma}\label{lemma:disjoint-partition}
    Let $G$ be a principal ample groupoid with compact unit space $\Go$. 
    For every compact open symmetric set $D\subseteq G$ and every $N \geq 1$, there exists a finite partition of compact open sets $\Go = C_1 \sqcup \cdots \sqcup C_k$ such that
    \begin{align*}
        B_D(u,N) \cap B_D (v,N) = \emptyset 
    \end{align*}
    for all $u \neq v \in C_i$ and $1 \leq i \leq k$. 
\end{lemma}

\begin{proof}
    As $G$ is \'etale, $\Go$ is open in $G$, and so $D^{2N} \setminus \Go$ is compact.
    Then the image $\mathcal{R} = (r,s)(D^{2N} \setminus \Go) \subseteq \Go \times \Go$ is compact, and since $\Go \times \Go$ is Hausdorff, we deduce that $\mathcal{R}$ is closed. Furthermore, since $G$ is principal, the image misses the diagonal of $\Go \times \Go$. 
    
    Fix $x\in \Go$. Since $(x,x)\notin\mathcal R$ and $\mathcal R$ is closed, there are open neighbourhoods $V_x,W_x$ of $x$ such that $(V_x\times W_x)\cap \mathcal R=\emptyset$. Using ampleness of $G$, we  find a compact open neighbourhood $O_x \subseteq V_x \cap W_x$ and observe
    \begin{align*}
        (O_x \times O_x ) \cap \mathcal{R} = \emptyset.
    \end{align*}
    Since $\Go$ is compact, we obtain a finite subcover $O_1, \ldots , O_k$. From this we may construct a disjoint compact open cover by setting
    \begin{align*}
 C_1=O_1,\qquad
 C_i=O_i\setminus\bigcup_{j<i}O_j
 \qquad \text{for $2\leq i\leq k$},
\end{align*}
Note that since  $C_i\subseteq O_i$, we have $(C_i\times C_i)\cap\mathcal R=\emptyset $. 

We claim that for $u \neq v \in C_i$, $1 \leq i \leq k$, we have $B_D(u,N) \cap B_D(v,N) = \emptyset$, which will finish the proof. Suppose for the sake of contradiction  that $u\neq v \in C_i$ and $z \in B_D(u,N) \cap B_D(v,N)$. 
By the definition of these balls,
there are elements $a,b\in D^N$ with
\begin{align*}
 s(a)=u,\qquad s(b)=v,\qquad r(a)=r(b)=z.
\end{align*}
The product $b^{-1}a$ then has source $u$ and range $v$, and by symmetry of $D$ we have $b^{-1}a \in D^{2N}$. Since $u \neq v$, $b^{-1}a \in D^{2N} \setminus \Go$. But then
\begin{align*}
 (v,u)=(r(b^{-1}a),s(b^{-1}a))
 \in\mathcal R\cap(C_i\times C_i),
\end{align*}
which contradicts $\mathcal R\cap(C_i\times C_i) = \emptyset$ from above. Thus $B_D(u,N) \cap B_D(v,N) = \emptyset$. 
\end{proof}

The following is our variant of the ball carving construction relative to a family of radius choice functions. Note that the choice functions do not need to satisfy any compatibility condition.
The flexibility of moving to a family of radius choice functions allows us to explicitly control the growth of the set $A$ relative to that of $Y$.

\begin{lemma}\label{lemma:separator-existence} Let $G$ be a principal ample groupoid with compact unit space. 
    Let further $N\geq 1$, let $A\subseteq G^{(0)}$ be compact open, and suppose that for every compact open $W\subseteq A$ a locally constant function
$
t_W:W\longrightarrow\{0,\ldots,N-1\}
$
is given. 

\noindent Then there exists a compact open $Y\subseteq A$ such that every $D$-component of $A\setminus Y$ is contained in a ball of radius $N-1$. In particular, $\left\langle D\cap G|_{A\setminus Y}\right\rangle\subseteq D^{2N}$.

\noindent Moreover, if $$\eta:=\sup\left\{\frac{\vert W\cap
\bigl(B_D(u,t_W(u)+1)\setminus B_D(u,t_W(u))\bigr)\vert}{\vert W\cap B_D(u,t_W(u))\vert} \colon W\subseteq A \text{ compact open}, u\in W\right\}$$
is finite, then $Y$ can be chosen to satisfy
$
\gamma_Y(R)\leq\eta\,\gamma_A(R+2N).
$
\end{lemma}

\begin{proof}
    If $A = \emptyset$, we use $Y = \emptyset$ and we are done. We may therefore assume $A \neq \emptyset$. Using Lemma \ref{lemma:disjoint-partition} we find a compact open partition $\Go = C_1 \sqcup \cdots \sqcup C_k$ such that for each $i$, $u \neq v \in C_i$ implies $        B_D(u,N) \cap B_D (v,N) = \emptyset$. 

    Set $W_0 = A$. For $i \geq 1$, suppose that $W_{i-1}$ has been defined and is compact open. Set $S_i = C_i \cap W_{i-1}$,
    and define
    \begin{align*}
        P_{i,u} &:= W_{i-1}\cap B_D(u, t_{W_{i-1}}(u)) \\
        Q_{i,u}&:= W_{i-1} \cap (B_D(u,t_{W_{i-1}}(u)+1) \setminus B_D(u,t_{W_{i-1}}(u))).
    \end{align*}
    Then by construction
    \begin{align}\label{eq:P-Q-props}
        u \in P_{i,u}, \qquad \vert Q_{i,u}\vert \leq \eta \vert P_{i,u}\vert , \qquad P_{i,u} \sqcup Q_{i,u} \subseteq B_D(u,N). 
    \end{align}

    For $0 \leq t < N$, put $E_{i,t} := \{u \in S_i \mid t_{W_{i-1}}(u) = t\}$. Since $t_{W_{i-1}}$ is locally constant, the sets $E_{i,t}$ are compact open. For a general compact open set $E \subseteq \Go$, we write
    \begin{align*}
        \mathcal{B}_n (E) := r(D^n \cap s^{-1}(E)) = \bigcup_{u \in E} B_D(u,n),
    \end{align*}
    which is again compact open by Lemma \ref{lemma:locally-constant-counts}. We deduce that the sets
    \begin{align*}
        P_i &:= \bigcup_{u \in S_i} P_{i,u} = W_{i-1} \cap \bigcup_{t=0}^{N-1} \mathcal{B}_t(E_{i,t} ), \text{ and } \\
        T_i &:= \bigcup_{u \in S_i} (P_{i,u} \sqcup Q_{i,u} ) = W_{i-1} \cap \bigcup_{t=0}^{N-1} \mathcal{B}_{t+1} (E_{i,t} )
    \end{align*}
    are compact open. Now observe that for $u\neq v \in S_i \subseteq C_i$ we must have $Q_{i,u} \cap P_{i,v} = \emptyset$ since $Q_{i,u} \subseteq B_D(u,N)$ and $P_{i,v} \subseteq B_D(v,N)$, and the latter sets are disjoint by Lemma \ref{lemma:disjoint-partition}. Together with $Q_{i,u} \cap P_{i,u} = \emptyset$ we deduce that
    \begin{align*}
        Q_i := \bigcup_{u \in S_i} Q_{i,u} = T_i \setminus P_i. 
    \end{align*}
    It follows that both $Q_i$ and $W_i := W_{i-1} \setminus T_i$ are compact open.

    The sets under consideration are pairwise disjoint across different stages of the induction, because the later stages take place in the remainder left after removing earlier pieces. Moreover, $W_i \cap C_i = (W_{i-1} \setminus T_i) \cap C_i = (W_{i-1} \cap C_i) \setminus T_i = S_i \setminus T_i = \emptyset$  because $S_i \subseteq P_i \subseteq T_i$. The remainders decrease and the $C_i$ cover $\Go$, so we deduce $W_k = \emptyset$. 
    Hence
    \begin{align}\label{eq:def-Y}
        Y := \bigcup_{i=1}^k Q_i \qquad \text{ satisfies } \qquad A\setminus Y = \bigsqcup_{i=1}^k \bigsqcup_{u\in S_i} P_{i,u},
    \end{align}
    and both $Y$ and $A \setminus Y$ are compact open. 

    We show that distinct pieces in the decomposition of $A \setminus Y$ can not be joined by an element of $D$. First, suppose $a \in P_{i,u}$ and $b \in P_{i,v}$ are joined by an element of $D$, and that $u \neq v$. Then $b \in B_D(u, t_{W_{i-1}}(u) +1) \subseteq B_D(u,N)$ while also $b \in B_D(v,N)$, which contradicts the disjointness of the balls provided by Lemma \ref{lemma:disjoint-partition}. If instead $a \in P_{i,u}$ and $b \in P_{j,v}$ are joined by an element of $D$ while $i < j$, then $b \in W_{j-1} \subseteq W_i$, but also
    \begin{align*}
        b \in W_{i-1} \cap B_D(u, t_{W_{i-1}}(u) +1 ) = P_{i,u} \sqcup Q_{i,u} \subseteq T_i,
    \end{align*}
    which contradicts $W_i \cap T_i = \emptyset$. We conclude that every $D$-component of $A\setminus Y$ is contained in some piece $P_{i,u} \subseteq B_D(u, N-1)$, which verifies the first claim of the lemma. The statement $  \langle D\cap G|_{A\setminus Y}\rangle \subseteq D^{2N}$ follows from  Lemma \ref{lemma:loc-to-glob-D-boundedness}. 

    We end by verifying the final claim. Fix $x \in \Go$ and $R \in \N_0$, and set
    \begin{align*}
        J := \{ (i,u) \mid 1\leq i \leq k, \, u \in S_i, \, Q_{i,u} \cap B_D(x,R) \neq \emptyset \}. 
    \end{align*}
    The $Q_{i,u}$ are pairwise disjoint, and $B_D(x,R)$ is finite, so we deduce that $J$ is finite. If $(i,u) \in J$, choose some $z \in Q_{i,u} \cap B_D(x,R)$. Then every point of $P_{i,u}$ can be reached from $x$ through a path via $z$ and $u$, which by construction then has length at most $R + (t_{W_{i-1}}(u) +1) + t_{W_{i-1}}(u) \leq R + 2N -1$. Then
    \begin{align*}
        \bigsqcup_{(i,u) \in J} P_{i,u} \subseteq A \cap B_D(x, R + 2N). 
    \end{align*}
    By using \eqref{eq:P-Q-props} we then calculate
    \begin{align*}
        \vert Y \cap B_D(x,R) \vert &= \sum_{(i,u) \in J} \vert Q_{i,u} \cap B_D(x,R) \vert  \leq \sum_{(i,u) \in J} \vert Q_{i,u} \vert \\
        &\leq \eta \sum_{(i,u) \in J} \vert P_{i,u} \vert \leq \eta \vert A \cap B_D(x, R + 2N) \vert. 
    \end{align*}
    This holds for each $x \in \Go$, and so taking the supremum over all $x \in \Go$ finishes the proof.  
\end{proof}

Our goal is now to produce a family of functions $(t_W)_W$ as in the hypothesis of the previous lemma such that $\eta$ is small. 
The form of $\eta$ of course immediately leads us to consider conditions under which the balls themselves form F\o lner sets. Hence the following subexponential growth type criterion should not come as a surprise:

\begin{lemma}\label{lemma:boundary-estimate}
    Let $N\in \mathbb N$ and $\varepsilon>0$ such that
    \begin{align}\label{eq:epsilon-N-condition}
        \gamma_{\Go}(N) < (1+ \varepsilon) ^N.
    \end{align}
    Then for every compact open $W\subseteq G^{(0)}$ the function $t_W:W\to \{0,\ldots, N-1\}$ given by
    $$t_W(u)=\min\{ t\colon  \vert W \cap (B_D (u, t+1) \setminus B_D(u,t)) \vert \leq \varepsilon \vert W \cap B_D(u,t) \vert\}$$
    is well defined and continuous.
\end{lemma}

\begin{proof}

    Set $f(t) = \vert W \cap B_D(u,t) \vert$, and observe that $f(0) = 1$ and $f(N) \leq \gamma_{\Go}(N)$. Note that for any $W \subseteq \Go$ and any $u \in W$, there exists some $0 \leq t < N$ which satisfies
  \begin{align}\label{eq:boundary-bound}
       \vert W \cap (B_D (u, t+1) \setminus B_D(u,t)) \vert \leq \varepsilon \vert W \cap B_D(u,t) \vert. 
    \end{align}
    Indeed, if the inequality \eqref{eq:boundary-bound} failed for each $0 \leq t < N$, then we would immediately get a contradiction to \eqref{eq:epsilon-N-condition}. The claim that the smallest such $t$ is locally constant then follows from Lemma \ref{lemma:locally-constant-counts}. 
\end{proof}

We are now ready to prove our main result:

\begin{proof}[Proof of Theorem~\ref{introthm:A}] We will prove the result under the additional assumption that $\Go$ is compact. As explained in section \ref{sec: prelims}, the general case follows by passing to the Aleksandrov groupoid $G^+$.

Fix a compact open and symmetric set $G^{(0)}\subseteq D\subseteq G$. Since $G$ has local polynomial growth of degree at most $s\geq 0$ we can find a constant 
    $C_D \geq 1$ such that
    \begin{align*}
        \sup_{x \in \Go} |D^n \cap G_x |\leq C_D (1+n)^s \quad \text{ for all } n \in \N_0
    \end{align*}
    Set $m = \lfloor s \rfloor +1$, and for $n\geq 2$ set
    \begin{align*}
        a_n = \frac{\log (2C_D(1+n)^s)}{n}, \quad \text{and} \quad \varepsilon_n = e^{a_n} -1.
    \end{align*}
     Now note that $\varepsilon_n = \mathcal{O}(\log(n)/n)$, and so
    \begin{align*}
        \varepsilon_n^m \gamma_{\Go} (2mn) \leq C_D \varepsilon_n^m (1 +2mn)^s = \mathcal{O}(n^{s-m} (\log(n))^m) \longrightarrow 0.
    \end{align*}
    We may therefore choose an $N\geq 1$ so large that 
    \begin{equation}\label{eq:controlled growth}
        \varepsilon_N^m \gamma_{\Go} (2mN) < 1.
    \end{equation}   
  By construction, $(1 + \varepsilon_N)^N = 2C_D (1+N)^s > \gamma_{\Go}(N)$, and so Lemma \ref{lemma:boundary-estimate} provides a family of locally constant functions $t_W:W\to \{0,\ldots, N-1\}$ such that
    $$\frac{\vert W\cap (B_D(u,t_W(u)+1)\setminus B_D(u,t_W(u)))\vert}{\vert W\cap B_D(u,t_W(u))\vert}\leq \varepsilon_N$$
    for every compact open $W\subseteq G^{(0)}$ and $u\in W$.
    Now set $A_0=\Go$ and inductively define a sequence of compact open sets $A_0\supseteq A_1\supseteq A_2 \supseteq \ldots$ by applying Lemma \ref{lemma:separator-existence} to $A_j$ to find $A_{j+1}\subseteq A_j$ such that
    \begin{enumerate}
        \item  $ \langle D\cap G|_{A_j\setminus A_{j+1}}\rangle \subseteq D^{2N}$, and 
        \item $\gamma_{A_{j+1}}(R) \leq \varepsilon_N \cdot \gamma_{A_j}(R+2N)$.
    \end{enumerate}
    We claim that $A_m=\emptyset$ so that setting $U_j:=A_j\setminus A_{j+1}$ for $j=0,\ldots, m-1$ yields the desired open cover of $G^{(0)}$ required by the definition of dynamic asymptotic dimension by item (1). To see that $A_m$ is indeed empty, repeatedly apply the estimate (2) to obtain 
    $$
    \gamma_{A_m}(0) \leq \varepsilon_N^m \cdot \gamma_{\Go}(2mN)\stackrel{\eqref{eq:controlled growth}}{<}1.
    $$ 
    Since $\gamma_{A_m}$ only takes non-negative integer values we conclude that $A_m=\emptyset$ as desired.
\end{proof}

Finally, we turn to the consequences concerning the dynamic asymptotic dimension growth. As the statement is slightly more technical, we repeat it for the convenience of the reader:
\begin{corollary}
    Let $G$ be a $\sigma$-compact, principal, ample groupoid with local polynomial growth and $f:[0,\infty)\to [0,\infty)$ an increasing unbounded function. Then there exists a continuous coarse length function $\ell$ on $G$ such that $\dad_{G,\ell}(t)\leq f(t)$. In particular, $\ell$ can be chosen so that $(G,\ell)$ has slow dynamic asymptotic dimension growth.
\end{corollary}
\begin{proof}
    We may assume that $\Go$ is compact by passing to $G^+$, as this preserves $\sigma$-compactness, principality, ampleness, and local polynomial growth. A continuous coarse length on $G^+$ restricts to one on $G$, and the dimension functions agree by our convention above. 
    Since $G$ is $\sigma$-compact we can find an increasing sequence
    $G^{(0)}\subseteq K_1\subseteq K_2\subseteq \ldots$ of compact open symmetric subsets exhausting all of $G$.
    Next, use that $G$ has local polynomial growth to find corresponding sequences of constants $C_n\geq 1$ and $s_n\geq 0$ such that for every $k\in \mathbb N$:
    $$\sup_{x \in \Go} \vert K_n^k\cap G_x\vert \leq C_n(1+k)^{s_n}.$$
    Let $G_n$ denote the subgroupoid of $G$ generated by $K_n$. Then $G_n$ has local polynomial growth of degree at most $s_n$. Indeed, given a compact open symmetric set $D\subseteq G_n$, there exists $q\in \mathbb N$ such that $D\subseteq K_n^q$. It follows that
    $$\sup_{x \in \Go} \vert D^k\cap (G_n)_x\vert \leq \sup_{x \in \Go} \vert K_n^{qk}\cap G_x\vert \leq C_n(1+qk)^{s_n}\leq C_nq^{s_n}(1+k)^{s_n}.$$
    It follows from Theorem~\ref{introthm:A} that $\dad(G_n)\leq \lfloor s_n\rfloor$. 
    Since $f$ is non-decreasing and unbounded, we can find a strictly increasing sequence $(a_n)_n$ of positive integers such that $f(a_n)\geq \lfloor s_n \rfloor$.
    Define a function $\ell:G\to [0,\infty)$ by
    $$\ell(g)=\min\left\{\sum_{j=1}^m a_{n_j} \middle| m\geq 1,\ g=g_1\cdots g_m,\ g_j\in K_{n_j}\right\}$$ for all $g\in G\setminus G^{(0)}$ and $\ell (g)=0$ for $g\in G^{(0)}$. It is routine to check that $\ell$ is a continuous coarse length function on $G$. 

    To complete the proof we show that $\dad_{G,\ell}(t)\leq f(t)$.  
    If $0 \leq t < a_1$, then $B_\ell(t) = \Go$ and $\dad_{G, \ell}(t) = 0 \leq f(t)$. For $t \geq a_1$, there is a unique $n \in \N$ such that $a_n \leq t < a_{n+1}$. 
    Suppose $g\in B_\ell(t)$ is a non-unit element. Let $g=g_1\cdots g_m$ be such that $\sum a_{n_j}$ realises $\ell(g)$.
    Then $\sum a_{n_j} = \ell(g)\leq t<a_{n+1}$. Thus, since $(a_n)_n$ is strictly increasing we must have $n_j\leq n$ and hence $g_j\in K_{n_j}\subseteq K_n$ for all $j=1,\ldots, m$.
    It follows that $B_\ell(t)\subseteq \langle K_n\rangle = G_n$. Since a cover witnessing finite dynamic asymptotic dimension at scale $B_\ell(t)$ in $G_n$ also witnesses finite dynamic asymptotic dimension at scale $B_\ell(t)$ in $G$ we get
    $$\dad_{G,\ell}(t)\leq \dad(G_n)\leq \lfloor s_n \rfloor\leq f(a_n)\leq f(t),$$
    as desired. 
\end{proof}

\section{Examples}
\label{sec:examples}
\subsection{Coarse groupoids and AF groupoids}
Let $(X,d)$ be a uniformly locally finite extended metric space satisfying $|B_d(x,R)| \leq C(1+ R)^s$ uniformly in $x \in X$. Its coarse groupoid $G(X,d)$ is principal and ample with unit space $\beta X$, where $\beta X$ is the Stone--\v{C}ech compactification of $X$ (see \cite{STY02} for more details on this construction). Moreover $\dad(G(X,d)) = \operatorname{asdim}(X,d)$ by \cite[Theorem 6.4]{GWY17}. Given a compact open symmetric $D \subseteq G(X,d)$, choose $R$ such that $D$ lies in the closure of the controlled set $E_R = \{(y,x) \mid d(x,y) \leq R \}$. For $x \in X$ we then have $|D^n \cap G(X,d)_x | \leq C(1 + nR)^s$. By Lemma \ref{lemma:locally-constant-counts} the function $\omega \mapsto |D^n \cap G(X,d)_\omega |$ is locally constant on $\beta X$, so density extends this bound to all $\omega \in \beta X$. Thus $G(X,d)$ has local polynomial growth of degree at most $s$, and Theorem \ref{introthm:A} gives $\operatorname{asdim}(X,d) \leq \lfloor s \rfloor$. For graphs with their extended graph metrics, this recovers the polynomial exponent bound in \cite{P23}. 

AF groupoids have local polynomial growth of degree zero. Indeed, by writing $G = \bigcup_j K_j$ as an increasing union of principal clopen subgroupoids with $K_j^{(0)} = \Go$ and $K_j \setminus \Go$ compact, every compact $D$ lies in some $K_j$, see e.g. \cite[Section 3]{AOP-JNCG}. Consequently $D^n \subseteq K_j$ and $\sup_{x \in \Go} |D^n \cap G_x | \leq \sup_{x \in \Go} |K_j \cap G_x| < \infty$ uniformly in $n$. Theorem \ref{introthm:A} therefore recovers $\dad(G) = 0$. 
\subsection{Partial action groupoids}
We refer the reader to \cite{AbadiePartialActions04} and \cite[Section 4.1]{AustadBonicke26} for material on partial action groupoids. 

Let $\theta=(X_g,\theta_g)_{g\in\Gamma}$ be a partial action of a discrete group $\Gamma$ on a locally compact totally disconnected Hausdorff space $X$ 
with open domains. The partial action $\theta$ is \emph{free} if $\theta_g(x)=x$ implies $g=e$ whenever defined.

Write $X\rtimes_\theta\Gamma$ for the associated partial transformation groupoid, whose elements are $(\theta_g(x),g,x)$ with $x\in X_{g^{-1}}$, see \cite[Section 4.1]{AustadBonicke26} and \cite{AbadiePartialActions04}. There is a natural identification of $(X\rtimes_\theta\Gamma)^{(0)} = X$. 

For a discrete group, local polynomial growth of degree at most $s$ means that for every finite symmetric set $F \subseteq \Gamma$ there is a constant $C_F\geq 1$ such that $|F^n| \leq C_F(1+n)^s$ for all $n \in \N_0$. 
\begin{corollary}\label{cor:partial-actions-dad}
Suppose $\Gamma$ is a discrete group with local polynomial growth of degree at most $s\geq 0$ and that $\theta=(X_g,\theta_g)_{g\in\Gamma}$ is a partial action on a locally compact totally disconnected Hausdorff space $X$ with open domains. If $\theta$ is free, then 
$$\dad (X\rtimes_\theta \Gamma)\leq \lfloor s \rfloor.
$$
\end{corollary}

\begin{proof}
    The groupoid $X \rtimes_\theta \Gamma$ is ample because $X$ is locally compact and totally disconnected, and it is principal because $\theta$ is free. The map
    \begin{align*}
        c\colon X \rtimes_\theta \Gamma \to \Gamma, \qquad c(\theta_g(x), g,x) = g,
    \end{align*}
    is a continuous groupoid homomorphism which is injective on each source fibre. If $D \subseteq X \rtimes_\theta \Gamma$ is compact open and symmetric, $F = c(D) \cup \{e\}$ is finite and symmetric and $c(D^n) \subseteq F^n$. Consequently,
    \begin{align*}
        | D^n \cap (X \rtimes_\theta \Gamma)_x | \leq |F^n | \leq C_F (1+n)^s
    \end{align*}
    for $x \in X$. Theorem \ref{introthm:A} now applies. 
\end{proof}

For \emph{global} actions of finitely generated groups the result above can also be obtained in a bit of a roundabout way by a combination of existing results in the literature: if a finitely generated group $\Gamma$ has polynomial growth of degree at most $d$, then by Gromov's Theorem it is virtually nilpotent. Combining \cite[Corollary~1.10]{Bartels17} and \cite[Theorem~4.7]{GWY17} we conclude that the dynamic asymptotic dimension of any free action of $\Gamma$ on a compact, metrizable, totally disconnected Hausdorff space is finite. Thus, \cite[Corollary~3.4]{Pilgrim25} shows $\dad(X\rtimes \Gamma)=\operatorname{asdim}(\Gamma)$ and the latter is at most $d$ by \cite{P23}. Thus, even though our result is not new in this situation, we believe it does provide a much more direct route.

\subsection{Point set groupoids}
In the interest of brevity, we refer the reader to \cite[Section 5.2]{AOP-JNCG} or \cite[Section 4.2]{AustadBonicke26} for material on point set groupoids. 

Let $H$ be a second-countable locally compact group.
 For a nonempty separated point set $\Lambda\subseteq H$, let $\Omega(\Lambda)$ be its Chabauty--Fell hull and put $\Omega_0(\Lambda)=\{P\in\Omega(\Lambda) \mid e\in P\}$.
As in \cite[Section 5.2]{AOP-JNCG}, its point set groupoid is
\begin{align*}
 G(\Lambda)&=\{(h,P)\in H\times\Omega_0(\Lambda):h^{-1}\in P\},\\
 s(h,P)&=P,\qquad r(h,P)=hP,
\end{align*}
with the topology and operations inherited from the translation groupoid. This is an \'etale groupoid.
We say that $\Lambda$ has \emph{finite local complexity} if for every compact $K \subseteq H$ the collection $\{(x^{-1} \Lambda) \cap K \mid x \in \Lambda \}$ is finite. 
We say $P \in \Omega(\Lambda)$ has \emph{trivial translation stabilizer} if $hP=P$ implies $h=e$. 

\begin{corollary}\label{cor:point-sets-dad}
Let $H$ be a second-countable locally compact group with left Haar measure $m$ and local polynomial growth of degree at most \(s\) in the sense that for every compact subset $K\subseteq H$, there is a constant $C_K>0$ such that
\[
m(K^n)\leq C_K(1+n)^s
\text{ for all }n\in\mathbb N.
\]

Suppose $\Lambda\subseteq H$ has finite local complexity and every nonempty
$P\in\Omega(\Lambda)$ has trivial translation stabilizer. 
Then $G(\Lambda)$ is principal and ample, and
\begin{align*}
 \operatorname{dad}(G(\Lambda))\leq\lfloor s\rfloor.
\end{align*}
\end{corollary}

\begin{proof}
The groupoid is Hausdorff and \'etale with compact unit space
$\Omega_0(\Lambda)$; see \cite[Section~5.2]{AOP-JNCG}. 
Finite local complexity makes this unit space zero-dimensional by the first part of the proof of \cite[Proposition~4.22]{AustadBonicke26}, which uses no relative
density. 
Thus $G(\Lambda)$ is ample. The stabilizer hypothesis makes it principal. Finally, we show that $G(\Lambda)$ has local polynomial growth of degree at most $s$. 
Let $D\subseteq G(\Lambda)$ be a compact open symmetric subset and let $K\subseteq H$ be the image of $D$ under the canonical projection map. Then $K$ is compact.
Choose a relatively compact open identity neighbourhood $V$ such that the translates $pV$, $p \in P$, are pairwise disjoint for every $P \in \Omega(\Lambda)$. Let $M\subseteq H$ be a compact set such that $K\cup V\subseteq M$ and $e\in M$.
Since the disjoint translates $pV$ for $p\in P\cap K^n$ are contained in $K^nV\subseteq M^{n+1}$, we have
$$\vert P\cap K^n\vert m(V)=\sum_{p\in P\cap K^n} m(pV)=m\left(\bigsqcup_{p\in P\cap K^n} pV\right)\leq m(M^{n+1}).$$
Using that $H$ has local polynomial growth we can find $C_M>0$ such that
$m(M^n)\leq C_M(1+n)^s$ for all $n\in \mathbb N$.
Thus, for every $P\in\Omega_0(\Lambda)$ we have
\begin{align*}
 | D^n \cap G(\Lambda)_P | \leq  | P\cap K^n | \leq \frac{m(M^{n+1})}{m(V)} \leq \frac{C_M}{m(V)} (2+n)^s\leq \frac{C_M 2^s}{m(V)}(1+n)^s.
\end{align*}
Thus, $G(\Lambda)$ has local polynomial growth of degree at most $s$ and Theorem \ref{introthm:A} applies. 
\end{proof}

\begin{remark}
    Note in particular that unlike in the comparison results of \cite{AustadBonicke26}, we are not requiring $\Lambda$ be relatively dense.
\end{remark}

\vspace{1cm} \noindent {\bf Use of AI assistance:} 
This project was initiated in February 2026 during a visit of the first author to Newcastle University. During this visit the authors developed a potential proof strategy for Theorem \ref{introthm:A} and established some preliminary results. 
At a later stage, the authors had conversations with ChatGPT-6.0 Astra to flesh out the details of the strategy developed in February. The technical details in the present proof are based on these conversations. As the arguments produced by AI were difficult to read and obscured many of the underlying ideas, the authors completely reorganised and rewrote the proof. Finally, AI tools were used for proofreading the manuscript.

\printbibliography

\end{document}